\documentclass{article}

\usepackage{graphicx}

\usepackage{amsfonts}
\usepackage{amsmath}

\newtheorem{theorem}{Theorem}
\newtheorem{acknowledgement}{Acknowledgement}

\newtheorem{corollary}{Corollary}[section]

\newtheorem{proposition}{Proposition}[section]

\newenvironment{proof}[1][Proof]{\textbf{#1.} }{\ \rule{0.5em}{0.5em}}

\begin{document}

\title{Shot noise distributions and selfdecomposability }
\author{Aleksander M. Iksanov and Zbigniew J.
Jurek\footnote{On leave from the University of Wroc\l aw, Wroc\l
aw, Poland.} \\
Cybernetics Faculty, Kiev T. Shevchenko National \\
University, Ukraine \\
Department of Mathematics, Wayne State University, \ \ \\
Detroit, MI 48202, USA \\
iksan@unicyb.kiev.ua and zjjurek@math.wayne.edu}
\maketitle

\begin{abstract}
Stationary (limiting) distributions of shot noise processes, with
exponential response functions, form a large subclass of positive
selfdecomposable distributions that we illustrate by many
examples. These shot noise distributions are described among
selfdecomposable ones via the regular variation at zero of their
distribution functions. However, slow variation at the origin of
(an absolutely continuous) distribution function is incompatible
with selfdecomposability  and this is shown in three examples.

\end{abstract}

\section{Introduction.}

Let us recall that by \emph{a general univariate shot noise process} we mean
a stochastic process of the form
\begin{equation}
X(t)=\sum_{\tau _{i}\leq t}R_{i}(t-\tau _{i}),\text{ }t\in R,
\text{ or }t\geq 0,
\end{equation}
\newline
where $\tau _{i}$'s are the points of a renewal process (or of a
point process) and $\{R_{i}(t),t\geq 0\},i=1,2,...$ is a sequence
of independent identically distributed (in short: i.i.d.)
measurable stochastic processes that are also independent of the
renewal (point) process $\tau_{i}$'s. The stochastic processes
$\{R_{i}(t),t\geq 0\},$ $i=1,2,...$ are called the \emph{response
processes}.

Shot noise processes have been used to model a lot of diverse
phenomena. Papers by Hsing \& Teugels (1989) and Vervaat (1979)
provide a variety of applications. For more recent results see,
Doney \& O'Brien (1991), Bondesson (1992) Samorodnitsky (1998),
 Kl\"{u}ppelberg et al. (2001).

Typically one discusses particular cases of the response
processes. A very common class  is the following one :
\begin{equation*}
R_{i}(t)=\xi _{i}h(t),\text{ }t\geq 0,
\end{equation*}
\newline
where $\{\xi _{i}\}$ is a sequence of independent and identically
distributed random variables and $h(t),t\geq 0,$ is a non-random
measurable function. This allows to separate a response amplitude
$\xi $ from the dynamics of a response given by a function $h$.
Thus we have
\begin{equation}
X(t)=\sum_{\tau _{i}\leq t}\xi _{i}h(t-\tau _{i}),\text{ }t\in R\text{ \ or }%
t\geq 0.
\end{equation}
Furthermore we will assume that in the process (2),
$\{\tau_{i},i\geq 1\}$ is a Poisson flow, i.e., arrival times in a
Poisson process.

Our main aim here is to study the stationary distributions of some
shot noise processes. In Section 2, we restrict our attention to
distributions of stationary versions of (2) with the exponential
response function $h(u)=\exp (-\omega u)$, where $\omega>0$ is a
fixed parameter. These still constitute a  quite large subclass of
the class of all positive selfdecomposable distributions, what is
illustrated by many examples. Shot noise distributions, among the
selfdecomposable ones are described by regular variation at zero
of their distribution functions (Theorem 1). In Section 3, we show
that slowly varying at zero (and absolutely continuous)
distribution functions are not selfdecomposable (Theorem 2). This
is illustrated by some examples. In Section 4, we present a
general result on existence of the shot-noise transform , although
their fixed points will be studied in a separate paper. Finally we
would like to stress that our primary interest is in probabilistic
(distributional) theoretical questions rather than in modelling a
realistic phenomena.

To simplify some terminology and notations, we will use of the
following abbreviations:

\emph{''d.f.''- distribution function, ''p.d.f.''- probability
density function, ''LT''-Laplace transform, ''i.i.d.''-
independent identically distributed, ''r.v.''- random variable,
''ID''- infinitely divisible, ''SD''- selfdecomposable or
selfdecomposability, ''$\overset{d}{=}$''- equality in
distribution, ''SN''- shot noise.}

\section{Exponential response and selfdecomposability.}

\strut The non-stationary SN process $X(t)$ with exponential impulse
response function is given by the equality
\begin{equation}
X(t):=X(0)e^{-\omega t}+\int_{0}^{t}e^{-\omega (t-u)}dA_{\xi }(u)\overset{d}{%
=}X(0)e^{-\omega t}+\int_{0}^{t}e^{-\omega u}dA_{\xi }(u),
\end{equation}
(the equality in distribution holds for each fixed t) where the
L\'{e}vy process $A_{\xi}(t)=\sum_{k=1}^{N_{\lambda }(t)}\xi _{k}$
is the compound Poisson process and $N_{\lambda }(t)$ is a Poisson
process with arrival times $\{\tau _{i}\},i\geq 1$ and intensity
$0<\lambda <\infty $. Furthermore $\xi $, $ \xi _{1},\xi _{2},...$
are \emph{positive} i.i.d. r.v., independent also of $N_{\lambda
}(t)$;  $0<\omega <\infty $ and $X(0)$ are given non-random
constants. Note that in the sequel we may drop the index $\xi $
(jump variable) from the compound Poisson process $A_{\xi }(t)$.

It is known that the limit in (3) exists, as $t\rightarrow \infty
$, if and only if $\mathbb{E}\ln (1+A_{\xi}(1))<\infty $. And this
is equivalent to $\mathbb{E}\ln (1+\xi)<\infty$; cf. Jurek \&
Vervaat (1983) or Jurek \& Mason (1993), p. 122 for Banach space
case. Furthermore, the limit of (3), as $t\rightarrow \infty $, is
the unique stationary distribution of the SN process; cf. Jurek \&
Mason (1993), Proposition 3.7.10, p. 161, where there is even more
general case discussed. Also cf. Sato \& Yamazato (1984).

Below $F(x)$ denotes the SN distribution of $X(\infty )$ in (3).
As it was mentioned above $F$ can be used to define a stationary
process
\begin{equation}
Y(t)=Y(0)e^{-\omega t}+\sum_{i=1}^{N_{\lambda }(t)}\xi _{i}e^{-\omega
(t-\tau _{i})},
\end{equation}
where $P\{Y(0)\leq x\}=F(x)$. However, this fact will not be
discussed or used in this paper.]

Recall that a probability distribution of a random variable $\eta
$ is called \emph{selfdecomposable} or \emph{L\'{e}vy class $L$
distribution} if, for every $0<c<1,$ there exists a rv $\eta_{c}$,
independent of $\eta$, such that the equality
\begin{equation*}
\eta \overset{d}{=} c \eta +\eta _{c },
\end{equation*}
holds true. A useful characterization of class $L$ (or $SD$)
distributions, crucial for this paper, is the following
\emph{random integral representation}:

\emph{the distribution of r.v. $ \eta $ is SD iff there exists a
unique, in distribution, L\'{e}vy process $Z$ with $E\ln (1+\left|
Z(1)\right| )<\infty $ \ and
\begin{equation}
\eta \overset{d}{=}\int_{0}^{\infty }e^{-u}dZ(u).
\end{equation}}
To the process $Z$ we refer to as the \emph{background driving
L\'{e}vy process}, (in short: BDLP) of $\eta$. Cf. Jurek \&
Vervaat (1983) or Jurek \& Mason, Chapter 3, or Sato \& Yamazato
(1984) or Wolfe (1982). The terminology of background driving
L\'{e}vy process was introduced in Jurek (1996).

A relation between positive SD r.v. $\eta $ with LT $\Phi $ and
its BDLP $Z(1)$ with LT $\Psi $, is the following one:
\begin{equation}
\log \Phi (s)=\int_{0}^{s}\log \Psi (r)\frac{dr}{r},\ \ \mbox{and}
\ \ \log\Psi(s)= \frac{s\Phi^{\prime}(s)}{\Phi(s)}.
\end{equation}
cf. Jurek \& Mason (1993) or Proposition 3 in Jurek (2001).
Further, let us  note that stationary (limiting) distributions of
SN processes are SD and their BDLP's are the compound Poisson
processes. This is because (4) can be rewritten as (5), with
$Z(t)$ being compound Poisson processes $A_{\xi}(t), t \ge 0$.

Before presenting next result let us recall that a function $L$ is
called \emph{slowly varying at 0} ( in Karamata sense), if for
each $\lambda> 0$, $\lim_{s \to 0}\frac{L(\lambda s)}{L(s)} = 1$.
In a case, the limit is of the form $\lambda^{\rho}$ we say that
$L$ is \emph{regularly varying with index $\rho$}. For the theory
of such functions we refer to Bingham, Goldie and Teugels (1987);
in short: BGT(1987).
\begin{theorem}
If $F(x)$ is the limiting distribution function of the SN process
(3) then it is regularly varying at 0 with index $\lambda
/\omega$, i.e.,
\begin{equation*}
F(x)\sim x^{\lambda /\omega}L(x),x\rightarrow 0\text{,}
\end{equation*}
with $L(x)$ being slowly varying at zero. Conversely, if a
positive $SD$ r.v. $\eta$ has the distribution function $F$
regularly varying at zero, with positive index, then it is $SN$
distribution.
\end{theorem}
\emph{Proof}. From the formula (3) we see that $X(\infty)$ is $SD$
r.v. with the compound Poisson process $A_{\xi}(t),t \ge 0$, as
its BDLP. Consequently, from Remark 3.6.9 (4) in Jurek \& Mason
(1993), p. 126, we get its LT
\begin{equation}
\Phi(s)=E\exp(-sX(\infty))=\exp
\{-(\lambda/\omega)\int_{0}^{s}(1-E\exp(-u\xi)du/u\}.
\end{equation}
Cf. also Jurek \& Vervaat (1983) for Banach space valued SD r.v.
or Jurek (1996),(1997),(2001). Now let us note that (7) can be
rewritten as follows
\begin{multline*}
\Phi(s)= s^{- \lambda/\omega}\exp[\lambda/\omega (\log s
-\int_{0}^{s}(1-E\exp(-u\xi))du/u)]= \\
s^{-\lambda/\omega}\Big(\exp[-\lambda/\omega\int_{0}^{1}(1-E\exp(-u\xi))du/u]\Big)
\exp[\lambda/\omega\int_{1}^{s}E(\exp(-u\xi))du/u].
\end{multline*}
\newline
By comparing the above with the Representation Theorem (Thm 1.3.1
in BGT (1987), p. 12) one observes that the LT $\Phi(s)$ of
$X(\infty )$ is regularly varying at infinity with index
$(-\lambda /\omega )$. Consequently, by Karamata's Tauberian
Theorem ( Thm 1.7.1, in BGT(1987) on p. 38), one gets that the
corresponding d.f. F is regularly varying at zero with the index
$\lambda/\omega$, which proves the first part of the theorem.

Conversely, suppose that $\eta >0$ is SD r.v., its LT $\Phi (s)$
is given by (6) and that the d.f. $P\{\eta \leq x\}$ is regularly
varying at zero with index $\rho > 0$, while its BDLP Z(t) is not
a compound Poisson process. Thus Z(1) has an infinite spectral
measure N and
\begin{equation*}
\underset{s\rightarrow
+\infty }{\lim }(-\ln \Psi (s))=\underset{s\rightarrow +\infty }{\lim }%
\int_{0}^{\infty }(1-e^{-sx})N(dx)=\infty ,
\end{equation*}
\newline
by the Lebesgue monotone convergence theorem . Note that the first
equality is just the representation of the LT for positive ID
distribution (the drift is zero because of regular variation
assumption). On the other hand, from (6), the Tauberian and the
Monotone Density Theorems (cf. BGT(1987) Thm 1.7.1 and 1.7.2) we
conclude
\begin{equation*}
\underset{s\rightarrow +\infty }{\lim }\ln \Psi (s)=\underset{%
s\rightarrow +\infty }{\lim }\dfrac{s\Phi^{\prime}(s)}{\Phi(s)}=
-\rho,
\end{equation*}
which is a contradiction and thus completes the proof.

[The above result, in Ukrainian, is given in Iksanov(2001) but
without a detailed proof.]

\begin{corollary} A positive SD rv with support $(0,\infty)$ and the
L\'evy spectral measure $M(dx)=k(x)/x dx$ is SN distributions if
and only if $k(0+)<\infty.$
\end{corollary}
\emph{Proof.} From Vervaat(1979) we know that $X(\infty)$ has
L\'evy spectral measure of the form
\begin{equation*}
M(dx)=(\lambda /\omega )P\{\xi>x\}dx/x.
\end{equation*}
This and Theorem 1 gives the proof.
\newline
\textbf{Remark.} Also in Vervaat (1979) we have that
\begin{equation*}
X(\infty )\overset{d}{=}R^{\omega /\lambda }(X(\infty )+\xi)
\end{equation*}
where rv $R$ has the the uniform distribution on $[0,1]$.
Consequently, Khintchine-Shepp criterion gives the unimodality of
$X(\infty)$ with mode at zero, provided   $\lambda /\omega \le 1$.
However, as we now know from Yamazato's (1978) result, all SD on
real line are unimodal.
\newline

\textbf{Examples of SN distributions.}

Before giving some
explicit examples of positive SN distributions, note that formula
(7) (or (6)) implies that
\begin{equation}
E[e^{-s\xi}]=1+(\omega/\lambda) s\{d\Big(\log
E[e^{-sx}dF(x)]\Big)/ds\},
\end{equation}
which allows the identification of \emph{the generating
distribution  of jumps} $\xi _{i}$, which appears in compound
Poisson process (BDLP) $A_{\xi}(.)$.
\newline

In the examples below we put $\rho =\lambda /\omega $ as a
parameter.

\textbf{1)} Gamma distributions or rv's  denoted by
$\gamma_{\rho,\beta}$, (where $\rho$  and $\beta$ are called the
{\emph{shape } and $\beta$ \emph{scale parameters}, respectively;
) are given via their p.d.f.
\begin{equation*}
f_{\rho ,\beta }(x)=\beta ^{\rho }/\Gamma (\rho )x^{\rho -1}\exp
(-\beta x),x>0,
\end{equation*}
are SN distributions  generated by the exponential  distributions
$P\{\xi \leq x\}=1-\exp(-\beta x)$. Cf. Jurek (1997), p. 97, for
more details.
\newline

\textbf{2)}  Distribution functions of the following form
\begin{equation}
F_{\rho ,\beta }(x)=1-\sum_{k=0}^{\infty }(-\beta )^{-k}x^{\rho k}/\Gamma
(1+\rho k),\text{ \ }x\geq 0,\text{ }0<\rho \leq 1,\text{ }\beta >0,
\end{equation}
are called \emph{positive Linnik distributions}. Their LT are
given by
\begin{equation}
\varphi _{\rho ,\beta }(s)=1/(1+\beta s^{\rho }).
\end{equation}
Using (8) we see that the positive Linnik distribution is
generated by itself. It is not hard to prove that these are the
only SN distributions with such a property (Iksanov (2001)).

The distributions with LT (10) are relatively well-known because
$\eta \overset{d}{=}Z(\varepsilon )$, where $\eta $ is the
positive Linnik distribution  with LT (10) and $Z(t)$ is a
positive strictly $\rho $-stable L\'evy process independent of
$\varepsilon$ that is an exponentially distributed r.v. with the
scale parameter $1/\beta $. (A \emph{'stochastic'} way of proving
SD of the above $\eta $ is presented in Jurek (2001), p. 244 .)
Also note that the explicit form (9) was given in Pillai (1990).
But let us add the following

\begin{corollary} For positive Linnik distributions
with $0<\rho <1$ we have that
\begin{equation*}
\underset{x\rightarrow \infty }{\lim }\dfrac{1-F_{\rho ,\beta
}(x)}{F_{\rho ,\beta }(1/x)}=\beta ^{2}\Gamma (1+\rho )/\Gamma
(1-\rho ).\text{\ }
\end{equation*}
\end{corollary}
\emph{Proof.} Apply Corollary 8.1.7 from BGT(1987), p. 334, to the
the upper tail of F, and the Karamata Tauberian Theorem 1.7.1 in
BGT(1987), p. 37, for the lower tail to conclude the proof.
\newline

\textbf{3)} \emph{A generalized positive Linnik distribution} is
given by the following LT
\begin{equation*}
\varphi _{\rho ,\rho _{1},\beta }(s)=1/(1+\beta s^{\rho
_{1}})^{\rho /\rho _{1}},0<\rho _{1}\leq 1,\text{ }\rho ,\beta >0.
\end{equation*}
Using (8) one notes it is generated by a positive Linnik
distribution with parameters $\rho _{1},\beta $.

Also, as before, a generalized positive  Linnik distributions
(with varying parameter $ \rho $) can be realized as a law of the
subordinated L\'evy process arising from a positive strictly
stable process and gamma process. ( Comp. Jurek (2001), p. 244).

\textbf{4)} When one randomizes  the shape parameter in the gamma
distribution $\gamma_{k+ \rho+1 , 1}$ , with $k$ distributed as
Poisson rv, one gets the so called \emph{Bessel distribution}; cf.
Feller (1966), p. 58. Explicitly, it is the distribution with
p.d.f. given by
\begin{equation*}
f_{\rho }(x)=\exp (-\rho -x)(x/\rho )^{(\rho -1)/2}I_{\rho
-1}(2\sqrt{\rho x} ),\text{ }x>0\text{, \ }\rho >0,
\end{equation*}
where $I_{\rho -1}(x)$ is the modified Bessel function with index
$\rho -1$. Bessel distribution is SN distribution when the
generating rv $\xi$ has gamma distribution, i.e., with d.f.
$P\{\xi \leq x\}=1-\exp (-x)(x+1)$. It follows from the fact that
the LT of Bessel distribution is of the form from the Remark below
with n=2.
\newline

\textbf{Remark.} 1). For BDLP of the form $A_{\gamma_{n,1}}(t)$ ,
in (3), one gets that
\begin{equation*}
E\exp (-sX(\infty ))=(1+s)^{-\rho }\prod_{i=1}^{n-1}\exp (-(\rho
/i)(1-(1+s)^{-i})),
\end{equation*}
i.,e., the law of $X(\infty ) $ is a convolution of
$\gamma_{\rho,1}$  with $n-1$ compound Poisson laws. The i-th of
them  has intensity $\rho /i$ and  gamma distribution ( with shape
parameter $i$ ) as the law of its jumps.

2) For $n=2$ the above reveals that Bessel distribution is a
convolution of SD distribution (i.e. gamma distribution) and its
BDLP evaluated at one . Two other such cases ( SD distribution
convoluted with its BDLP gives SD distribution) are given in Jurek
(2001), p. 248. Both are related to some integral functionals of a
Brownian motion .
\newline

In all SN distributions \textbf{1)-4)} we know explicitly the
generating (the jumps) distribution  $\xi$, while for those
\textbf{5)-7)} below we do not. For the examples \textbf{5)} and
\textbf{6)}, jumps $\xi$ have d.f. with completely monotone
derivatives ; cf. Bondesson (1992), p. 68 and p. 59 property (iv)
for the example \textbf{5)} and p. 60 for example \textbf{6)}.
\newline

\textbf{5)} \emph{Burr distributions or a generalized beta
distributions of second type have d.f.} are given by
\begin{equation*}
F_{\rho ,\beta _{1},\beta _{2}}(x)=1-(\beta _{1}/(x^{\rho }+\beta
_{1}))^{\beta _{2}}\text{, }x\geq 0,\text{ }\beta _{1}>0,\beta _{2}>\rho >0.
\end{equation*}

These include Pareto when $\rho =1$ and $F$-distribution for
properly specified values of the parameters. From Theorem 1 it
follows immediately  that those are indeed SN distributions.
\newline

\textbf{6)} \emph{ Weibull distribution} is defined by d.f.
\begin{equation*}
F_{\rho ,\beta }(x)=1-\exp (-\beta x^{\rho }),\text{ }x\geq
0,0<\rho \leq 1, \text{ }\beta >0;
\end{equation*}
\newline

\textbf{Remark.} Except for the example \textbf{4)}, all the above
distributions have completely monotone p.d.f.'s, whenever $\rho
\leq 1.$ Goldie (1967) and Steutel (1970) proved that all
non-negative r.v.'s whose d.f.'s have completely monotone
derivatives are ID (infinitely divisible). However, they are not
necessarily SD even after replacing ''derivative'' by ''p.d.f.''.
As a counterexample take a distribution with p.d.f.
\begin{equation*}
g_{a,b}(x)=\exp (-ax)(1-\exp (-x))^{b-1}/B(a,b),
\end{equation*}
where $b\leq 1$ and $2a+b\geq 1$ ; c.f. Bondesson (1992), p. 143.
Similar examples are discussed in the following section.
\newline

\textbf{Remark.} Note that if the distribution with completely
monotone p.d.f. is SD then it can realized as for $X(\infty )$,
with $\rho \leq 1$. On the other hand, even for $\rho \leq 1$ the
limiting laws of (3) do not necessarily have completely monotone
p.d.f.'s . These is so because these p.d.f.'s can have
discontinuous derivatives, as is in case of degenerated at 1 jumps
$\xi_i$ .
\newline

\textbf{7).} Absolute value of Cauchy rv  leads to
\emph{half-Cauchy distribution with p.d.f.}
\begin{equation*}
f_{1}(x)=2(\pi (1+x^{2}))^{-1}\text{, }x>0\text{,}
\end{equation*}
that by Theorem 1 gives SN distribution. Comp. Di\'{e}dhiou(1998),
where the selfdecomposability was originally proved via an
analytic argument.
\newline

Finally let us remark that all positive selfdecomposable
distributions are weak limits of sequences of laws of $X(\infty
)$, i.e., shot noise distributions.
\newline
\section{Positive non-selfdecomposable distributions.}

In this section we provide some tools for checking the
non-selfdecomposabilty of a positive and absolutely continuous rv.
We begin with the following.

\begin{theorem}
If a positive rv $\eta$ has the distribution function $F(x)$
slowly varying at zero then $F(x)$ is not selfdecomposable.
\end{theorem}
\emph{Proof.} Recall that SD distributions  have LT of the from
$\exp\int_{0}^{\infty}(e^{-sx}-1)\frac{k(x)}{x}dx$, with function
$k(x)$ being positive and non-increasing. From Theorem 1 and
Corollary 2.1., $k(0+)$ is finite iff the BDLP is compound
Poisson. On the other hand, $k(0+)=\infty$ implies that the
corresponding d.f. F is rapidly varying at zero ( by rapidly we
mean: regularly varying with index $\infty$). This completes the
proof.

\textbf{Examples of non-selfdecomposable distributions.}

\textbf{A).} For a fixed $\epsilon>0$, define a measure
\begin{equation*}
M_{\epsilon}(dx)= (x\log^{1+\epsilon}(x))^{-1}1_(x>e)(x)dx +
x^{-1}1_{(0,e]}(x)dx ,
\end{equation*}
which is a L\'evy spectral measure of SD d.f., say $F_{\epsilon}$;
( note that $xdM_{\epsilon}(x)/dx$ is non-increasing function).
Since the upper tail $M_{\epsilon}( [x, \infty))$ is a slowly
varying , then so is $1-F_{\epsilon}(x)$, by Thm 8.2.1 in
BGT(1987). Now let define rv $\eta_{\epsilon}$ whose tail is equal
to LT of SD rv with d.f. $F_{\epsilon}$, i.e.,
\begin{equation*}
P(\eta_{\epsilon} >x)=\int_{0+}^{\infty }\exp
(-xy)dF_{\epsilon}(y)\text{.}
\end{equation*}
Applying  Corollary 8.1.7 in BGT(1987) we see it is also slowly
varying function. So, in view of the above Theorem 2 we conclude
that $\eta_{\epsilon} $ is not  SD.
\newline

\textbf{Remark.} Let us recall here, that p.d.f.'s which are
"proportional" to LT's of SD distributions are SD; cf. Bondesson
(1992), p. 28.
\newline

\textbf{B).} Before stating our next result  let us recall that
strictly stable processes "stopped" at independent positive SD rv
give again SD rv.( cf. Jurek (2001), Proposition 1, or Bondesson
(1992).) In other words, L\'evy processes on $R$ subordinated to a
strictly stable process by a SD non-decreasing L\'evy process
(subordinator) is SD. Sato (2000) proved analogous result for a
Brownian motion with drift as the subordinand. However, this is
not true for non-Gaussian stable processes with drift. From
Ramachandran (1997), Theorem 2.1, p. 302, we conclude that the
L\'evy process $Y_{1}(Y_{2}(t))$ is \emph{not} SD , where
$Y_{2}(t)$ is the gamma process and $Y_{1}(t)$ is the stable with
the characteristic function $E\exp(izY_{1}(t))= \exp t(i\mu
z-b\left| z\right| ^{\alpha })$, $1<\alpha <2$, $\mu $, $b>0$.
Next proposition provides similar cases.

\begin{proposition}
\strut If $X_{1}(t)$ and $X_{2}(t)$ are two independent
non-decreasing L\'evy SD processes , whose BDLP's are compound
Poisson processes, then the subordinated  L\'evy process
(subordinator) $ X_{3}(t)=X_{1}(X_{2}(t))$ is not SD.
\end{proposition}

\begin{proof}
\strut It suffices to verify the case of $t=1$. Put $\Phi
_{i}(s)=E\exp (-sX_{i}(1))$, $i=1,2,3$. These LT's are related as
follows: $\Phi _{3}(s)=\Phi _{2}(\ln \Phi _{1}(s))$. From Theorem
1 we infer that there exist $\alpha _{1}$,$\alpha _{2}>0$ such
that
$\Phi_{i}(s)\sim s^{-\alpha _{i}}L_{i}(s)$, $s\rightarrow \infty $, where $%
L_{i}(s)$ are slowly varying at infinity, $i=1,2$. Formula (7)
with $X_{1}(1)$ instead
of $X(\infty )$ and $\alpha _{1}$ instead of $\lambda /\mu $ implies $%
\underset{s\rightarrow \infty }{\lim }s\Phi _{1}^{^{\prime
}}(s)/\Phi _{1}(s)=-\alpha _{1}.$ Now we can use de
L'H\^{o}pital's rule to get
\begin{equation*}
\underset{s\rightarrow \infty }{\lim }\dfrac{-\ln \Phi
_{1}(s)}{\alpha _{1}\ln s}=\ \ \underset{s\rightarrow \infty
}{\lim }\dfrac{-s\Phi _{1}^{^{\prime }}(s)/\Phi _{1}(s)}{s(\alpha
_{1}/s)}=1.
\end{equation*}
Thus $\ln \Phi _{1}(s)$ is a slowly varying. Condition
$\underset{s\rightarrow \infty }{\lim }\ln \Phi _{1}(s)=-\infty $,
guarantees the slow variation of $\Phi _{3}(s)$ as the
superposition of regularly and slowly varying functions; cf. BGT
(1987) Prop.1.5.7 . This with Theorem 2 completes the proof.
\end{proof}
\newline

\textbf{C).} Consider the log-Cauchy law with p.d.f.
$f(x)=(\pi x(1+ \ln^{2}x))^{-1}$, $x>0$. Di\'{e}dhiou (1998) points out
without proof that this law is not SD. But this  immediately
follows from our Theorem 2.
\newline

\textbf{D).} Here is an example of a phenomena that
\emph{"stopping" SD L\'evy process at non-infinitely divisible
random time  may get SD distribution.} Let $S_{1/2}(t)$ be a
strictly stable subordinator with the index of stability $1/2$ and
let $s_{1/2,1}$, $s_{1/2,2}$ be independent copies of
$S_{1/2}(1)$. Recall that the $1/2$ - stable d.f. coincides with
d.f. of $N^{-2}$, where N is a normal rv. Then the r.v.
$s_{1/2,1}/s_{1/2,2}\overset{d}{=} S_{1/2}(s_{1/2,2}^{-1/2})$ has
p.d.f. $p(x)=(\pi \sqrt{x}(x+1))^{-1}$, $x>0$. This is p.d.f. of
the beta distribution of the second kind ( ratio of two
independent gamma rv ) and, hence, it is SD; cf. Bondesson (1992),
p. 59 . On the other hand, we will verify that the r.v.
$s_{1/2,2}^{-1/2}$ has LT $T_{1/2}(\lambda): = \sum_{k=0}^{\infty
}(-1)^{k}\lambda^{k}(\Gamma (1+\alpha k))^{-1}$ and hence is not
ID because the tail of its distribution tends to zero more rapidly
than it is possible for a non-normal i.d. distribution; see
Theorem 8.1.12 of BGT (1987). To this end, let $\eta _{1/2}$ be a
r.v. having positive Linnik distribution with parameters $1/2$,
$1$ and let $\varepsilon $ be an independent exponentially
distributed r.v. with scale $1$. Then
\begin{equation*}
\eta _{1/2}\overset{d}{=}S_{1/2}(\varepsilon )
\overset{d}{=}s_{1/2,2}\varepsilon ^{2} \ \ \mbox{or} \ \ \eta _{1/2}^{1/2}\overset{d}{=}%
\varepsilon /s_{1/2,2}^{-1/2}.
\end{equation*}
This is equivalent to the following
\begin{equation*}
T_{1/2}(x)=P\{\eta _{1/2}^{1/2}>x\}=\int_{0}^{\infty }\exp
(-xy)dP\{s_{1/2,2}^{-1/2}\leq x\}.
\end{equation*}
For other examples of similar kind cf. Kelker (1971) and Sato
(2000).
\newline

\section{A remark on a general shot noise transform.}
To study fixed points of shot noise transforms one needs to find the
corresponding domains. We will study the problem of fixed point in a
separate paper, while in this section we focus on existence of some random
integrals or infinite series of independent rv. \newline

Fix probability space $(\Omega $, $\mathcal{F}$, $\mathbb{P})$ and
let $ \mathcal{P}(\mathbb{R}^{d})=\mathcal{P}$ denotes the set of
all Borel probability measures on Euclidean space
$\mathbb{R}^{d}$, let $A_{\xi }(t)$ denotes the compound Poisson
process, with the jump rv $\xi $,and $0<\lambda <\infty $ being
the intensity parameter. Finally, let $h:(0,\infty )\rightarrow
\mathbb{R}^{d}$ be a fixed measurable function and $t>0$ fixed
parameter . Now let us define \emph{random integral mappings}
\begin{equation*}
\mathcal{I}_{h}^{(t)}:\mathcal{P}\rightarrow ID\ \ \mbox{and}\ \
\mathcal{I} _{h}^{(t)}(\mu _{\xi
}):=\mathcal{L}(\int_{0}^{t}h(s)dA_{\xi }(s))=\mathcal{L
}(\sum_{\tau _{i}\leq t}h(\tau _{i})\xi _{i})=\nu _{t}\in ID,
\end{equation*}
where $\mu _{\xi }=\mathcal{L}(\xi )$. As in the beginning of
Section 2, we refer to Jurek \& Vervaat (1983), Lemma 1.1 or Jurek
\& Mason (1993), Chapter 3 , to conclude that $\nu _{t}$ is $ID$
and that its \emph{standard } L\'{e}vy-Khintchine representation
has the following ''triplet'' : a shift $ a_{h}^{(t)}$, a zero
Gaussian part and the L\'{e}vy spectral measure $ M_{h}^{(t)}$,
and
\begin{multline*}
a_{h}^{(t)}=\int_{0}^{t}\int_{\mathbb{R}^{d}}\frac{h(s)x}{1+h^{2}(x)||x||^{2}
}\mu _{\xi }(dx)ds=\int_{0}^{t}\mathbb{E}\Big[\frac{h(s)\xi
}{1+h^{2}(x)\xi
^{2}}\Big]ds, \\
M_{h}^{(t)}(B):=\lambda \int_{0}^{t}T_{h(s)}\mu _{\xi
}(B)ds=\lambda \int_{0}^{t}\mu _{\xi }(\frac{1}{h(s)}B)ds=\lambda
\int_{0}^{t}\mathcal{L} (h(s)\xi )(B)ds,
\end{multline*}
where $B$ is any Borel subset of $\mathbb{R}^{d}-\{0\}$ and $t$ is fixed.

Note, however, because of the compound Poisson input $A_{\xi
}(t)$, the stochastic process
\[
t \to \int_{0}^{t}h(s)dA_{\xi}(s), \ \  t\ge 0,
\]
is with independent increments and has paths with bounded
variation. Consequently, by Gihman and Skorohod (1975), Chapter
IV, Theorem 8, p. 279,  $\nu_t$  has "finite variation"
L\'evy-Khintchine representation $(0,0,M_{h}^{(t)})$, with the
integrability condition $\int_{\mathbb{R}^{d}}(\left\| x\right\|
\wedge 1)M_{h}^{(t)}(dx)<\infty$, for each $t>0.$

\begin{proposition}
In order that $X(\infty ):=\lim_{t\rightarrow \infty
}\int_{0}^{t}h(s)dA_{\xi }(s)=\sum_{k=1}^{\infty }h(\tau _{k})\xi _{k}$
exists in distribution (or in probability or a.s.) it is necessary and
sufficient that
\begin{eqnarray}
\int_{\mathbb{R}^{d}}(1\wedge ||y||)M_{h}^{(\infty )}(dy)
&=&\int_{0}^{\infty }\int_{\mathbb{R}^{d}}(1\wedge |h(s)|||y||)\mu
_{\xi
}(dy)ds=  \notag \\
&=&\int_{0}^{\infty }\mathbb{E}(1\wedge |h(s)|||\xi ||)ds<\infty .
\end{eqnarray}
\end{proposition}
\textbf{Proof.} Denoting by $<.,.>$ a scalar product in
$\mathbb{R}^d$, we note that
\begin{multline}
\mathbb{E}\Big[\exp i<z,\int_{0}^{t}h(s)dA_{\xi }(s)>\Big]=\exp
\int_{0}^{t}\log \phi _{A_{\xi }(1)}(h(s)z)ds=  \notag \\
=\exp \int_{0}^{t}\lambda \int_{\mathbb{R}^{d}}(e^{i<z,h(s)x>}-1)\mu _{\xi
}(dx)ds=\exp \int_{\mathbb{R}^{d}}(e^{i<z,x>}-1)M_{h}^{(t)}(dx)\text{.}
\end{multline}
is a ch.f. of an ID distribution, without Gaussian part and a
shift, and whose L\'{e} vy spectral measure $M_{h}^{(t)}$
satisfies the condition
\begin{equation}
\int_{\mathbb{R}^{d}}(\left\| x\right\| \wedge 1)M_{h}^{(t)}(dx)<\infty
\text{, for each }t>0\text{.}
\end{equation}
Since, for each fixed $t>0$, $M_{h}^{(t)}=(\mathbb{P}\times
Leb)\circ H^{-1}$, where $H:\Omega \times (0,\infty )\rightarrow
\mathbb{R}^{d}$ is given by $H(w,s):=\xi (w)h(s)$ therefore the
above is equivalent to
\begin{equation*}
\int_{0}^{t}\mathbb{E}(1\wedge |h(s)|||\xi ||)ds<\infty \text{,
for each }t>0 \text{.}
\end{equation*}
L\'{e}vy spectral measures $M_{h}^{(t)}(.)$ are increasing in t.
The limit exists, as $t\to \infty$, iff (12) is satisfied  for
$t=\infty$. Moreover,we infer that the limit is also without
Gaussian part; cf. Araujo-Gine (1980), Corollary 5.3 in Chapter 2
or Theorem 4.7 in Chapter 3. Finally to see that all the three
modes of convergence are equivalent use either Jurek-Vervaat
(1983) or a theorem on convergence of infinite series of
independent summands. This completes the proof.

\begin{corollary}
An improper random integral $\int_{0}^{\infty }h(s)dA_{\xi }(ds)$
or an infinite series $\sum_{k=1}^{\infty }h(\tau _{k})\xi _{k}$
exists a.s. ( or in distribution or in probability ) if and only
if $\int_{0}^{\infty } \mathbb{E}(1\wedge |h(s)|||\xi ||)ds<\infty
$. Furthermore, its characteristic function is of the form $\exp
\int_{0}^{\infty }\lambda \int_{\mathbb{R}
^{d}}(e^{i<z,h(s)x>}-1)\mu _{\xi }(dx)ds$ .
\end{corollary}

\textbf{Remark} Our proof above is valid also for Hilbert and
Banach space valued r.v.'s.

\textbf{Examples.} Below we use the condition from Corollary 4.1 but the
explicit calculations are not given.

1) For $h(s)= e^{-s}$ one has that $\sum_{k=1}^{\infty}
e^{-\tau_k}\xi_k$ exists iff $\mathbb{E}[\log^{+}||\xi||]< \infty$
iff $\mathbb{E} [\log(1+||\xi||)] < \infty$.

2) For $h(s)=s^{-\alpha}, \alpha >1$, we $\sum_{k=1}^{\infty}
\tau_k^{-\alpha}\xi_k$ exists iff $\mathbb{E}[||\xi||^{1/\alpha}]< \infty$.

3) For $h(s)= 1_{(0<s<1)}(s)$, the infinite series
$\sum_{k=1}^{\infty} 1_{(0<s<1)}(\tau_k) \xi_k$ converges always.
\newline

\emph{The univariate general shot noise transform} $\mathbb{T}_h$ is defined
as follows
\begin{equation}
\mathbb{T}_h(\mathcal{L}(\xi))=\mathcal{L}(\sum_{i=1}^{\infty}
\xi_{i}h(\tau_{i}))=\mathcal{L}(\int_{0}^{\infty}h(x)dA_{\xi}(dx)),
\end{equation}
for $\xi$ such that series or integral converge. We say that $\xi$
is \emph{ a fixed point} of the shot noise transform
$\mathbb{T}_h$ whenever $\mathcal{ L}(\xi) =
\mathbb{T}_h(\mathcal{L}(\xi))$. Those fixed points are, of course
ID, and their existence and properties will be studied in a
separate paper.

\begin{acknowledgement}
A part of this research was initiated during  A. Iksanov's visit
to Wroc\l aw. He gratefully acknowledges the hospitality and
financial support received from the University of Wroc\l aw,
Wroc\l aw, Poland.
\end{acknowledgement}

\end{document}